\documentclass{article}
\usepackage{amsmath,amsthm,amssymb}
\usepackage[T1]{fontenc}

                         \newtheorem{Th}{Theorem}
                         \newtheorem{Prop}{Proposition}
                         \newtheorem{Lemma}{Lemma}
                          \newtheorem{Cor}{Corollary}
\theoremstyle{remark}
\newtheorem{Remark}{Remark}
\newtheorem*{Example}{Example}
\newtheorem{Definition}{Definition}
\newtheorem*{prooftwgl}{Proof of Theorem~\ref{twgl}}

\newcommand{\cb}{\mathcal{B}}
\newcommand{\ca}{\mathcal{A}}

\newcommand{\cs}{\mathcal{S}}

\newcommand{\vep}{\varepsilon}

\newcommand{\un}{\underline}

\newcommand{\rat}{\operatorname{R}}

\newcommand{\beq}{\begin{equation}}
\newcommand{\eeq}{\end{equation}}
\newcommand{\R}{{\mathbb{R}}}
\newcommand{\T}{{\mathbb{T}}}
\newcommand{\cir}{{\mathbb{S}^1}}
\newcommand{\C}{{\mathbb{C}}}
\newcommand{\Z}{{\mathbb{Z}}}
\newcommand{\N}{{\mathbb{N}}}
\newcommand{\Q}{{\mathbb{Q}}}
\newcommand{\xbm}{(X,\mathcal{B},\mu)}

\newcommand{\wne}{\operatorname{Int}}
\newcommand{\var}{\operatorname{Var}}
\begin{document}

\renewcommand{\thefootnote}{\fnsymbol{footnote}}

\title{Mild mixing property for special flows under piecewise constant functions}
\author{K.\ Fr\k{a}czek, M.\
Lema\'nczyk \\
{\small Faculty of Mathematics and Computer Science} \\
{\small Nicolaus Copernicus University}\\
{\small Ul. Chopina 12/18, 87-100 Toru\'n, Poland}
\\ E.\ Lesigne\\
\small{Laboratoire de Mathématiques et Physique Théorique}\\
\small{Faculté des Sciences et Techniques}\\
\small{ Université François Rabelais de Tours}\\
\small{ parc de Grandmont, 37200 Tours, France} }

\maketitle
\renewcommand{\thefootnote}{}
 \footnote{The first two authors were partially
supported by KBN grant 1 P03A 03826}
\begin{abstract}
We give a condition on a piecewise constant roof function and an irrational
rotation by $\alpha$ on the circle to give rise to a  special flow having the
mild mixing property. Such flows will also satisfy Ratner's property. As a
consequence we obtain a class of mildly mixing singular flows on the two--torus
that arise from quasi-periodic Hamiltonians flows by velocity changes.
\end{abstract}

\section{Introduction}
A finite measure--preserving dynamical system is mildly mixing (see
\cite{Fu-We}) if its Cartesian product with an arbitrary ergodic conservative
(finite or infinite) measure-preserving dynamical system remains ergodic.  It
is an immediate observation that the (strong) mixing of a dynamical system
implies its mild mixing.

The purpose of this paper is to present a certain class of special
flows over circle irrational rotations which are mildly mixing but
not mixing. The first example of mildly mixing dynamical system
that is not mixing (using Gaussian processes) was constructed by
Schmidt in \cite{Sch}. Another famous example of a mildly mixing
but not mixing system is the well--known Chacon transformation
(here the mild mixing property follows directly from the minimal
self-joining property). Recently, the first two authors proved in
\cite{Fr-Le} that the class of special flows built from piecewise
absolutely continuous functions with a non--zero sum of jumps and
over a rotation by $\alpha$ with bounded partial quotients is
mildly mixing. The absence of mixing follows from a more general
result of  Kochergin  \cite{Ko} concerning special flows built
over any irrational rotation on the circle and under functions of
bounded variation.

Opposite to the situation in  \cite{Fr-Le}, in this paper we will
deal with the case where the sum of jumps of the roof function $f$
is zero. Let $Tx=x+\alpha$ be an irrational rotation on the
circle, where $\alpha$ has bounded partial quotients. Let
$f:\T\to\R$ be a positive piecewise absolutely continuous function
whose derivative is square integrable and let
$\{\xi_1,\ldots,\xi_p\}$ stand for the set of its discontinuity
points. Let $d_i=d(\xi_i)$ stand for the value of the jump at
$\xi_i$, $i=1,\ldots,p$, i.e.\
\[d(\xi)=\lim_{x\to\xi^-}f(x)-\lim_{x\to\xi^+}f(x).\]
Suppose that the sum of jumps of $f$, i.e.\ $S(f)=\sum_{i=1}^pd_i$, is equal to
zero. Then without loss of generality we can assume that $f$ is piecewise
constant. Indeed, this follows from the fact that every absolutely continuous
function whose derivative is square integrable is cohomologous (whenever
$\alpha$ has bounded partial quotients) to a constant function (the proof of
this fact involves the classical small divisor argument).

Let us define now two additional properties  that $f$ may satisfy:
\begin{description}
\item[(P1)] whenever a nontrivial sum $\sum_{i=1}^p n_id_{i}$ with
$n_i\in\Z$, $i=1,\ldots,p$ equals zero, there exist $1\leq
i<i'\leq p$ such that $n_i$ and $n_{i'}$ are nonzero and
$\xi_i-\xi_{i'}\in(\Q+\Q\alpha)\setminus(\Z+\Z\alpha)$;
\item[(P2)]
\[\int f(x)\,dx\notin\sum_{i=1}^p(\xi_i+\Q+\Q\alpha) d_i.\]
\end{description}

\begin{Example}
If $f=a+b\chi_{[0,\xi)}$, where $a,b>0$ and $a/b\notin\Q+\Q\alpha$ and
$\xi\in(\Q+\Q\alpha)\setminus(\Z+\Z\alpha)$, then $f$ satisfies (P1) and (P2).
\end{Example}

The main result of the paper, whose proof is postponed until
Section~\ref{secmm}, is the following condition for the mild mixing.

\begin{Th}\label{twgl}
If $\alpha$ has bounded partial quotients and $f:\T\to\R^+$ is a
piecewise constant function satisfying properties (P1) and (P2)
then  $T^f$ is mildly mixing.
\end{Th}

In the case where the sum of jumps is non--zero the proof of the
mild mixing property (see \cite{Fr-Le}) is based upon two tools: a
property that emulates the Ratner condition introduced in
\cite{Ra} (we will describe it in detail in
Section~\ref{sectionratner}); and the absence of partial rigidity.

In the case of zero sum of jumps only the first property will
survive, although Ratner's property will follow from different
arguments than those from \cite{Fr-Le}; the property (P1) will
play the main role. In contrast with the method from \cite{Fr-Le},
special flows considered here are partially rigid (special flows
under piecewise constant functions are partially rigid for an
arbitrary irrational $\alpha$; see Remark~\ref{partrig}) and we
can not apply Lemma~3 in \cite{Fr-Le}.

In Section~\ref{secint} we deal with a natural class of singular flows on the
two--torus that arise from quasi-periodic Hamiltonians flows by velocity
changes. We give some conditions, under which they have special representations
satisfying (P1) and (P2).

\section{Basic definitions and notation}
\subsection{Factors}
Assume that $\cs=(S_t)_{t\in\R}$ is an ergodic flow on a standard
probability space $\xbm$. By that we mean always a so called {\em
measurable flow}, i.e.\ we require that the map $\R\ni t\to\langle
f\circ S_t,g\rangle\in\C$ is continuous for each $f,g\in L^2\xbm$.
Let $(t_n)_{n\in\N}$ be a sequence of real numbers such that
$t_n\to+\infty$. We say that a flow $\cs$ on $\xbm$ is {\em rigid}
along $(t_n)$ if
\begin{equation}\label{defszt}
\mu(A\cap S_{-t_n}A)\to \mu(A)
\end{equation}
for every $A\in\mathcal{B}$. By a {\em factor} we will mean any
$\cs$--invariant sub--$\sigma$--algebra $\mathcal{A}$ of
$\mathcal{B}$. Then a factor $\mathcal{A}\subset\mathcal{B}$ of
$\cs$ is rigid along $(t_n)$ if the convergence (\ref{defszt})
holds for every $A\in\mathcal{A}$.

Recall that a flow is mildly mixing iff it has no non--trivial
rigid factor (see~\cite{Aar} p.84).

A flow $\cs$ on $\xbm$ is called {\em partially rigid} along
$(t_n)$ if there exists $0<u\leq 1$ such that
\[\liminf_{n\to\infty}\mu(A\cap S_{-t_n}A)\geq u\mu(A)\;\;\text{ for every }\;\;A\in\mathcal{B}.\]

\begin{Lemma}\label{al5}
Let $\cs=(S_t)_{t\in\R}$ be an ergodic measure--preserving flow on
a probability standard Borel space $\xbm$. Suppose that $\cs$ is
not mildly mixing. Then there exists a nontrivial factor
$\mathcal{A}\subset\mathcal{B}$ of $\cs$ and $(t_n)_{n\in\N}$ with
$t_n\to+\infty$ such that
\[S_{t_n}\to E(\,\cdot\,|\mathcal{A}),\]  weakly in
$\mathcal{L}(L^2(X,\mu))$, i.e.\ for any $g,h\in L^2(X,\mu)$ we
have
\[\lim_{n\to\infty} \int_X h\circ S_{t_n}\cdot g \,d\mu=\int_X E(h|\mathcal{A}) \cdot g
\,d\mu.\]
\end{Lemma}
This lemma is a continuous time version of a result which can be
found in \cite{Le-Pa} for measure--preserving transformations. For
the completeness we outline the main arguments from \cite{Le-Pa}.
The proof is based on the theory of joinings of dynamical systems.
For the background of this theory we refer the reader to
\cite{Gl}. By a {\em self--joining} of a flow $\cs=(S_t)_{t\in\R}$
on $\xbm$ we mean any probability $(S_t\times
S_t)_{t\in\R}$--invariant measure on $(X\times
X,\mathcal{B}\otimes {\mathcal{B}})$ whose projections on $X$ are
equal to $\mu$. The set of self--joinings of $\cs$, denoted by
$J(\cs)$, considered  with the weak topology becomes a metrizable
compact semitopological semigroup for the operation of composition
of self--joinings (see \cite{Gl} Ch.6 for definitions). For every
$t\in\R$ by $\mu_{S_t}\in J(\cs)$ we define the graph joining
determined by $\mu_{S_t}(A\times B)=\mu(A\cap S_t^{-1}B)$ for
$A,B\in\mathcal{B}$. Suppose that $\mathcal{A}'$ is a non-trivial
rigid factor of $\cs$. Then we put
\[ I(\mathcal{A}')=\{\rho\in
J(\cs):\:
\rho|_{\mathcal{A}'\otimes\mathcal{A}'}=(\mu_{Id})_{\mathcal{A}'}\;\mbox{
and }\; \rho\in\overline{\{\mu_{S_t}:\:t\in\R\}}\}.
\]
Then, as noticed in \cite{Le-Pa}, $I(\mathcal{A}')$ is a
semitopological compact semigroup and hence it contains an
idempotent, which corresponds to a factor fulfilling the
requirements of Lemma~\ref{al5}.

Suppose that $\mathcal{A}\subset\mathcal{B}$ is a factor. Let us consider the
standard (quotient) probability Borel space
$(X/\mathcal{A},\mathcal{A},\mu_{\mathcal{A}})$, where $\mu_{\mathcal{A}}$ is
the image of $\mu$ via the factor map
$\pi:(X,\mathcal{B})\to(X/\mathcal{A},\mathcal{A})$. Denoting by
$\mathcal{P}(X,\mathcal{B})$ the set of probability measures on
$(X,\mathcal{B})$ let $\{\mu_{\un{x}}\in\mathcal{P}(X,\mathcal{B}):\un{x}\in
X/\mathcal{A}\}$ be the disintegration of $\mu$ over the factor
$\sigma$--algebra $\mathcal{A}$. Then (see \cite{Fur})
\[\mu_{\un{x}}(A)=E(\chi_A|\mathcal{A})(\un{x})\]
 and
\[\mu(A)=\int_{X/\mathcal{A}} \mu_{\un{x}}(A)\,d\mu_{\ca}(\un{x})\]
for every $A\in\mathcal{B}$ and for $\mu_{\mathcal{A}}$--a.e.\ $\un{x}\in
X/\mathcal{A}$. We say that $\mathcal{A}$ has {\em finite fibers} if there
exists a natural number $k$ such that $\mu_{\un{x}}$ is the uniform measure on
a $k$-element set for $\mu_{\mathcal{A}}$--a.e.\ $\un{x}\in X/\mathcal{A}$.

\subsection{Fraction expansion and discrepancy}
We denote by $\T$ the circle group $\R/\Z$ which we will
constantly identify with the interval $[0,1)$ with addition mod
$1$. Let $\lambda$ stand for Lebesgue measure on $\T$. For a real
number $t$ denote by $\{t\}$ its fractional part and by $\|t\|$
its distance to the nearest integer number. For an irrational
$\alpha\in\T$ denote by $(q_n)$ its sequence of denominators (see
e.g.\ \cite{Ch}), that is we have
\begin{equation}\label{ulla}
\frac{1}{2q_nq_{n+1}}<\left|\alpha-\frac{p_n}{q_n}\right|<\frac{1}{q_nq_{n+1}},
\end{equation}
where
\[\begin{array}{ccc}
q_0=1, & q_1=a_1, & q_{n+1}=a_{n+1}q_n+q_{n-1}\\
p_0=0, & p_1=1,  & p_{n+1}=a_{n+1}p_n+p_{n-1}
\end{array}\]
and $[0;a_1,a_2,\dots]$ stands for the continued fraction
expansion of $\alpha$.
\begin{Prop}\label{nierdon}
(Denjoy--Koksma inequality, see \cite{Her}) If $f:\T\to\R$ is a
function of bounded variation then
\[\left|\sum_{k=0}^{q_n-1}f(x+k\alpha)-q_n\int_{\T}f\,d\lambda\right|\leq\var f\]
for every $x\in \T$ and $n\in\N$.
\end{Prop}
 We say that $\alpha$ has {\em bounded
partial quotients} if the sequence $(a_n)$ is bounded. If
$C=\sup\{a_n:n\in\N\}+1$ then
\[\frac{1}{2Cq_n}<\frac{1}{2q_{n+1}}<\|q_n\alpha\|<\frac{1}{q_{n+1}}<\frac{1}{q_n}\]
for each $n\in\N$.

The following lemma is well--known.
\begin{Lemma}\label{al1}
Let $\alpha\in\T$ be irrational with bounded partial quotients.
Then there are constants $C_1,C_2>0$ such that for each $m\in\N$
the lengths of intervals $I$ in the partition of $\T$ arisen from
$0,-\alpha,\ldots,-(m-1)\alpha$ satisfy
$$
\frac{C_1}{m}\leq |I|\leq\frac{C_2}{m}.$$
\end{Lemma}

We will also need  the following simple fact.

\begin{Lemma}\label{al2}
Let $\alpha\in\T$ be irrational with bounded partial quotients and
let $\beta\in\Q+\Q\alpha$. Then there exists $c_1>0$ such that for
each $m\in\N$ the lengths of intervals $I$ in the partition of
$\T$ arisen from
$0,-\alpha,\ldots,-(m-1)\alpha,\beta,\beta-\alpha,\ldots,\beta-(m-1)\alpha$
is at least $c_1/m$.
\end{Lemma}

\begin{proof}
Assume that $\beta=(p_1+p_2\alpha)/q$, where $p_1,p_2\in\Z$ and
$q\in\N$. Suppose that $0\leq r,s<m$ and $\beta-s\alpha\neq
-r\alpha\text{ mod }1$. If $q(r-s)+p_2\neq 0$, then from
Lemma~\ref{al1} we have
\begin{eqnarray*}
q\left\|\left(\beta-s\alpha\right)+r\alpha\right\|&\geq
&\left\|q\left(\left(\beta-s\alpha\right)+r\alpha\right)\right\|=
\|(q(r-s)+p_2)\alpha\|\\
&\geq & \frac{C_1}{qm+|p_2|}\geq\frac{C_1}{(q+|p_2|)m},
\end{eqnarray*}
whence $\left\|(\beta-s\alpha)+r\alpha\right\|\geq
\frac{C_1}{q(q+|p_2|)m}$. If $q(r-s)+p_2=0$, then
\[\left\|\left(\beta-s\alpha\right)+r\alpha\right\|=\left\|\frac{p_1}{q}\right\|\geq\frac{1}{q}\geq\frac{1}{qm}.\]
\end{proof}

\section{Singular flows on two--torus}
\label{secint}

Let $H:\R^2\to\R$ be a $C^{\infty}$--quasi--periodic function,
i.e.\
\[H(x+m,y+n)=H(x,y)+m\alpha_1+n\alpha_2\]
for all $(x,y)\in\R^2$ and $m,n\in\Z$, and where
$\alpha=\alpha_1/\alpha_2$ is irrational. Points of the two-torus
are denoted as $\overline x$, with real local coordinates $x$ and
$y$.  Then $H$ determines a (quasi--periodic) Hamiltonian flow
$(h_t)_{t\in\R}$ on the torus associated with the following
differential equation
\[\frac{d\overline{x}}{dt}=X_H(\overline{x}),\text{ where
}X_H=\left(\frac{\partial H}{\partial y},-\frac{\partial
H}{\partial x}\right).\]

\begin{figure}[h]
\begin{center}
\vspace{5.2cm}\includegraphics{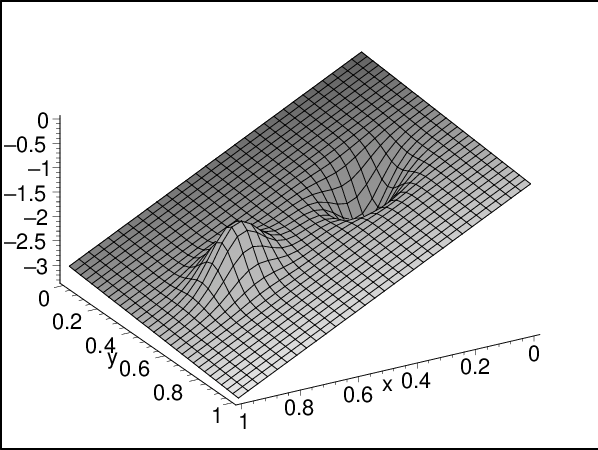}
\end{center}\caption{}
\label{wykham}
\end{figure}

Suppose that $H$ has critical points and $H$ is in the general
position, i.e.\ $H$ has no degenerate critical points and has all
critical values distinct. An example of the graph of such a
Hamiltonian is plotted in Figure~\ref{wykham}. In particular, each
critical point is either a non--degenerate saddle point or a
non--degenerate center. Moreover critical point repeats
periodically (with period 1 in each coordinate) but their critical
values are distinct.

\begin{figure}[h]
\begin{center}
\vspace{6cm} \includegraphics{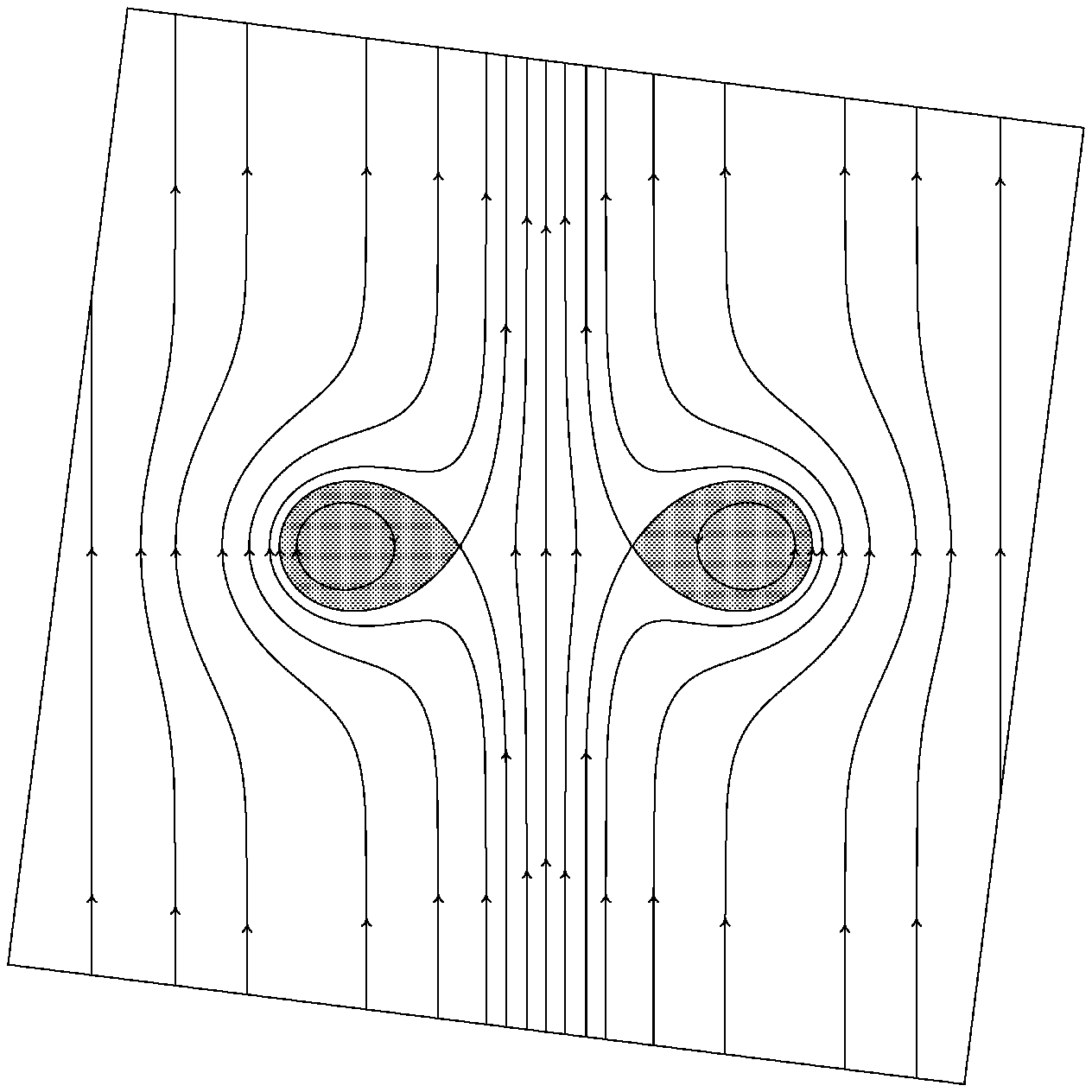}
\end{center}\caption{}
\label{twotraps}
\end{figure}

Then recall results by Arnold \cite{Arn}. Any connected component
of a level set of $H$ passing through a critical point is either
bounded (a point or a lemniscate--like curve) or it has the shape
of a folium of Descartes. In the unbounded case, the critical
value level set of $H$ separates the plane into two unbounded
components and a disk; the closure of the disk is called a {\em
trap}. A trap is homeomorphic to a closed disk and has a critical
point on the boundary, called the vertex of the trap (the same
terminology supplies when we pass to $\T^2$). The phase portrait
of the Hamiltonian flow on the torus determined by the
quasi--periodic Hamiltonian from Figure~\ref{wykham} is presented
in Figure~\ref{twotraps}. In this case the flow has two traps
(shaded areas). Then traps with distinct vertices are disjoint.
Therefore the phase space of $(h_t)_{t\in\R}$ decomposes into
traps filled with fixed points, separatrices and periodic orbits,
and an ergodic component $EC$ (white area in
Figure~\ref{twotraps}) of positive Lebesgue measure.

\subsection{From $(h_t)$ to singular flows}
Now will change velocity in the flow $(h_t)_{t\in\R}$. Let $B_i$,
$i=1,\ldots,p$ be all traps and $\overline{x}_i$, $i=1,\ldots,p$ their
vertices. Suppose $g:\T^2\to\R$ is a non--negative $C^{\infty}$--function which
is positive on the torus except of the points
$\{\overline{x}_1,\ldots,\overline{x}_p\}$. Let us consider the flow
$(\varphi_{t})_{t\in\R}$ on
$\T^2\setminus\{\overline{x}_1,\ldots,\overline{x}_p\}$ associated with the
following differential equation
\[
\frac{d\overline{x}}{dt}=X(\overline{x}), \text{ where
}X(\overline{x})=\frac{X_{H}(\overline{x})}{g(\overline{x})}.\]
Since the orbits of $(\varphi_{t})$ and $(h_t)$ are the same
(modulo fixed points of $(h_t)$), the phase space of
$(\varphi_{t})_{t\in\R}$ decomposes into traps filled with
critical points, separatrices and periodic orbits, and the ergodic
component $EC$ with positive Lebesgue measure.

Let us denote by $\omega=\omega_X\in {\mathcal{F}}^1(EC)$ the
$1$--form of class $C^{\infty}$ given by $\omega(Y)={\langle
X,Y\rangle}/{\langle X,X\rangle}$.

As it was proved in \cite{Fr-Le1}, $\int_{EC}d\omega_X$ exists,
and if $\alpha$ has partial bounded quotients and
$\int_{EC}d\omega_X\neq 0$ then the flow $(\varphi_t)_{t\in\R}$ on
$EC$ is weakly mixing and simple in the joining sense.

The aim of this paper is an exploration of the case
$\int_{EC}d\omega_X=0$. Let us consider the induced $1$--form
$i_{X_H}\nu$ on $\T^2$, where $\nu=dx\wedge dy$, i.e.\
\[i_{X_H}\nu(Y)=dx\wedge dy(X_H,Y).\]
Since the flow $(h_t)$ preserves the form $\nu$, $d(i_{X_H}\nu)=0$
(Liouville theorem). Let $\sigma_1,\sigma_2:[0,1]\to\T^2$ stand
for the pair of generators of $H_1(\T^2,\Z)$ given by
\[\sigma_1(t)=(t,0),\;\;\sigma_2(t)=(0,t)\text{ for }t\in[0,1].\]
Then
\begin{eqnarray*}
\int_{\sigma_1}i_{X_H}\nu&=&\int_0^1dx\wedge
dy(X_H(t,0),(1,0))dt\\&=&\int_0^1\frac{\partial}{\partial
x}H(t,0)dt=H(1,0)-H(0,0)=\alpha_1
\end{eqnarray*}
and similarly $\int_{\sigma_2}i_{X_H}\nu=\alpha_2$.

Fix $\overline{x}_0\in\T^2$. Then
$\xi_i=\int_{\overline{x}_0}^{\overline{x}_i}i_{X_H}\nu$ is well
defined up to $\Z\alpha_1+\Z\alpha_2$ for $i=1,\ldots,p$.
Moreover, for every  $i=1,\ldots,p$ put
\[d_i=-\int_{\partial B_i}\omega_X.\]
As it was shown in \cite{Fr-Le1} $d_i$ is well defined for
$i=1,\ldots,p$ and
\[\sum_{i=1}^p d_i=\int_{EC}d\omega_X=0.\]

\begin{Th}\label{twgladkie}
Suppose that $\alpha_1/\alpha_2$ is an irrational number with
bounded partial quotients. Assume that $\int_{EC}d\omega_X=0$. If
$X$ satisfies the following properties:
\begin{enumerate}
\item whenever a nontrivial sum $\sum_{i=1}^p n_id_{i}$ with
$n_i\in\Z$, $i=1,\ldots,p$ equals to zero,  there exist $1\leq
i<i'\leq p$ such that $n_i$ and $n_{i'}$ are nonzero and
$\xi_i-\xi_{i'}\in(\Q\alpha_1+\Q\alpha_2)\setminus(\Z\alpha_1+\Z\alpha_2)$;
\item \[\int_{EC}g\,d\nu\notin\sum_{i=1}^p(\xi_i+\Q\alpha_1+\Q\alpha_2) d_i,\]
\end{enumerate}
then the flow $(\varphi_t)$ is mildly mixing.
\end{Th}

\begin{Remark}
Notice  that $\xi_i-\xi_{i'}+\Z\alpha_1+\Z\alpha_2$ does not depend on the
choice of $\overline{x}_0$. Since $\sum_{i=1}^pd_i=0$,
$\sum_{i=1}^p(\xi_i+\Z\alpha_1+\Z\alpha_2) d_i$ does not depend on the choice
of $\overline{x}_0$ as well.
\end{Remark}

\subsection{Special representation of $(\varphi_t)$}
In this section we will describe a special representation of the
flow $(\varphi_t)$ and we will indicate why
Theorem~\ref{twgladkie} follows from Theorem~\ref{twgl}.

Assume that $T$ is an ergodic automorphism of a standard
probability space $\xbm$. If $f:X\to\R$ is a strictly positive
integrable function, then by $T^f=(T^f_t)_{t\in\R}$ we will mean
the special flow under $f$ (see e.g.\ \cite{Co-Fo-Si}, Chapter 11)
acting on $(X^f,{\cal B}^f,\mu^f)$, where $X^f=\{(x,s)\in X\times
\R:\:0\leq s<f(x)\}$ and ${\cal B}^f$ $(\mu^f)$ is the restriction
of ${\cal B}\otimes{\cal B}(\R)$ $(\mu\otimes\lambda_\R)$ to $X^f$
($\lambda_\R$ stands for Lebesgue measure on $\R$). Under the
action of the flow $T^f$ each point in $X^f$ moves vertically at
unit speed, and we identify the point $(x,f(x))$ with $(Tx,0)$.
More precisely, if $(x,s)\in X^f$ then
\[T^f_t(x,s)=(T^nx,s+t-f^{(n)}(x)),\]
where $n\in\Z$ is a unique number such that
\[f^{(n)}(x)\leq s+t<f^{(n+1)}(x)\]
with
\[f^{(m)}(x)=\left\{\begin{array}{ccc}
f(x)+f(Tx)+\ldots+f(T^{m-1}x) & \mbox{if} & m>0\\ 0 & \mbox{if} &
m=0\\ -\left(f(T^mx)+\ldots+f(T^{-1}x)\right)  & \mbox{if} & m<0.
\end{array}\right.\]

We can now describe a special representation of the flow
$(\varphi_t)$ which appeared in Section~\ref{secint}.  As it was
proved by Arnold in \cite{Arn}, there exists a closed
$C^{\infty}$--curve in $EC$ transversal to the orbits of
$(h_t)_{t\in\R}$, and which is homologous to $\sigma_2$. Fix a
point $\overline{x}_0$ on the transversal. Let
$\sigma:[0,\alpha_2]\to\T^2$ be the induced parameterization of
the transversal such that $\sigma(0)=\overline{x}_0$ and
\[\int_{\sigma(0)}^{\sigma(s)}i_{X_H}\nu=s\text{ for all }s\in[0,\alpha_2].\]
Moreover, the first--return map (Poincar\'{e} map) is determined everywhere on
the curve, except for the finite set $\{\sigma(\beta_i):i=1,\ldots,p\}$ of
points that are points of the last intersection of the incoming separatrices of
vertices of traps with the transversal curve ($\sigma(\beta_i)$ corresponds to
$\overline{x}_i$). In the induced parameterization, the first--return map is
the rotation by $\alpha_1$ mod $\alpha_2$. Moreover,
$\beta_i\in\xi_i+\Z\alpha_1+\Z\alpha_2$ for $i=1,\ldots,p$.

For the flow $(\varphi_t)_{t\in\R}$ on $EC$ we will consider the
same transversal. Then $\sigma$ is also the induced
parameterization for the vector field $X$ and the form
$g(x,y)\,dx\wedge dy$ which is preserved by $X$. As it was proved
in \cite{Fr-Le1}, the return time is a $C^{\infty}$--function  of
the parameters from $[0,\alpha_2]$ everywhere except of points
$\beta_i$, $i=1,\ldots,p$. Moreover, this function
$f:[0,\alpha_2]\to\R^+$ is piecewise absolutely continuous with
jumps at $\beta_i$ of size $d_i$ for $i=1,\ldots,p$  and the
derivative of $f$ is square integrable. Next notice that
\[\int_{EC}g(x,y)\,dx\,dy=\int_0^{\alpha_2}f(x)\,dx.\]
Indeed, let us consider the parameterization
\[\theta:\{(s,x):x\in[0,\alpha_2),0\leq s<f(x)\}\to EC\;\text{ given
by }\;\theta(s,x)=\varphi_s(\sigma(x)).\] Then
\[\frac{d}{ds}\theta(s,x)=X(\theta(s,x))\]
and
\[\int_{EC}g(x,y)\,dx\,dy=\int_0^{\alpha_2}\int_0^{f(x)}
g(\theta(s,x))|\det D\theta(s,x)|\,ds\,dx.\] Moreover
\begin{eqnarray*}
\det
D\theta(s,x)&=&\det\left[\frac{d}{ds}\theta(s,x)\;\frac{d}{dx}\theta(s,x)\right]=
\det\left[X(\theta(s,x))\;\frac{d}{dx}\theta(s,x)\right]\\
&=&\frac{1}{g(\theta(s,x))}DH(\theta(s,x))\frac{d}{dx}\theta(s,x)=\frac{1}{g(\theta(s,x))}\frac{d}{dx}(H\circ\theta)(s,x).
\end{eqnarray*}
Moreover $H(\theta(s,x))=H(\varphi_s(\sigma(x)))=H(\sigma(x))$ and
\[x=\int_{\sigma(0)}^{\sigma(x)}i_{X_H}\nu=\int_0^x\nu(X_H,\sigma'(t))\,dt=\int_0^x\frac{d}{dt}H\circ\sigma(t)\,dt=
H(\sigma(x))-H(\sigma(0)).\] Consequently,
\[\int_{EC}g\,d\nu=\int_0^{\alpha_2}\int_0^{f(x)}
\,ds\,dx=\int_0^{\alpha_2}f(x)\,dx.\]

Finally, rescaling $[0,\alpha_2]$ by the multiplication by $1/\alpha_2$ we
conclude that $(\varphi_t)$ is isomorphic to the special flow built over the
circle rotation $Tx=x+\alpha$, where $\alpha=\alpha_1/\alpha_2$ and under a
piecewise absolutely continuous function $f:\T\to\R$ whose derivative is square
integrable. Moreover, the sum of jumps $S(f)$ equals zero and if
$\Xi=\{\xi_1,\ldots,\xi_p\}$ is the set of discontinuity points of $f$ and
$d_i=d(\xi_i)$ stand for the value of the jump at $\xi_i$, $i=1,\ldots,p$, then
$f$ satisfies the properties (P1) and (P2). Without loss of generality we can
assume that $f$ is right-continuous and  $0\leq \xi_1<\ldots<\xi_p<1$.

\section{Properties of special flows under piecewise constant roof functions}
Let $Tx=x+\alpha$ be an irrational rotation and let $f:\T\to\R^+$
be a piecewise absolutely continuous and right--continuous
function such that $S(f)=0$. Let $\Xi=\{\xi_1,\ldots,\xi_p\}$
stand for the set of discontinuity points of $f$ and
$d_i=d(\xi_i)$ stand for the value of the jump at $\xi_i$,
$i=1,\ldots,p$.

\begin{Remark}
Recall that if $\alpha$ has bounded partial quotients and
$\zeta:\T\to\R$ is an absolutely continuous function with zero
mean such that $\zeta'\in L^2(\T,\lambda)$ then by the classical
small divisor argument $\zeta$ is a coboundary. It follows that
every piecewise absolutely function $f:\T\to\R$ whose derivative
is square integrable and such that $S(f)=0$ is cohomologous to a
piecewise constant function whose discontinuities and jumps agree
with discontinuities and jumps of $f$. We will henceforth assume
that $f$ is piecewise constant.
\end{Remark}

\begin{Remark}\label{nieciagl}
Notice that if $f$ satisfies the property (P1) then for every
$m\in\N$ the set $\{\xi_i-j\alpha:1\leq i\leq p,0\leq j<m\}$ is
the set of discontinuities of $f^{(m)}$. Indeed, if $f^{(m)}$ has
no jump at $\xi_i-j\alpha$ then there exist $1\leq
k_1<\ldots<k_u\leq p$ and $0\leq j_l<m$ for $l=1,\ldots,u$ such
that $\xi_{k_l}-j_l\alpha=\xi_i-j\alpha$ for $l=1,\ldots,u$ and
$\sum_{l=1}^u d_{k_l}=0$, which is impossible because
$\xi_{k_l}-\xi_{k_{l'}}\in\Z+\Z\alpha$ for all $1\leq l,l'\leq u$.
\end{Remark}

\begin{Remark}\label{remnien}
Let $a$ stand for the minimal value of $f$.  Notice that
\[\int f(x)\,dx\in a+\sum_{i=1}^pd_i\xi_i+\sum_{i=1}^p\Z d_i.\]
Indeed, since $\sum_{i=1}^pd_i=0$, we have
\[\int f(x)\,dx=f(\xi_p)-\sum_{i=1}^{p-1}\left(\sum_{j=1}^{i}d_j\right)(\xi_{i+1}-\xi_i)=a+(f(\xi_p)-a)+\sum_{i=1}^pd_i\xi_i.\]
Our claim follows from the fact that differences of all values of
$f$ belong to $\sum_{i=1}^p\Z d_i$. Using additionally (P2) it
follows that
\begin{equation}\label{wznien}
a\notin\sum_{i=1}^p\Q d_i.
\end{equation}
\end{Remark}

We will now show that (P1) together with (P2) imply the weak mixing of $T^f$.
Recall (see e.g.\ \cite{KaR}) that all we need to show for the weak mixing of
$T^f$ is that for every $r\in\R\setminus\{0\}$ the equation
\[\psi(Tx)/\psi(x)= e^{2\pi irf(x)}\]
has no measurable solution $\psi:\T\to\cir$. Let
$\sim\subset\Xi\times \Xi$ stand for the equivalence relation
given by $x\sim y$ iff $y-x\in\Z\alpha$. Since
$\sum_{i=1}^pd_i=0$, the property (P1) implies the existence of at
least two nontrivial equivalence classes.

\begin{Prop}(see Proposition~1.9~(iii) and Corollary~1.6 in \cite{Ga-Le-Li})\label{kawkoc}
The equation
\[\psi(Tx)/\psi(x)= e^{2\pi if(x)}\]
has a measurable solution $\psi:\T\to\cir$ iff
\[\sum_{\xi'\sim\xi}d(\xi')\in\Z\text{ for each }\xi\in\Xi\text{ and }\int f(x)\,dx\in \Z+\Z\alpha.\]
\end{Prop}

\begin{Cor}\label{slabemieszanie}
If $f$ satisfies (P1) and (P2) then the special flow $T^f$ is
weakly mixing.
\end{Cor}

\begin{proof}
Suppose contrary to our claim that $T^f$ is not weakly mixing. Then by
Proposition~\ref{kawkoc} there exists $r\in\R\setminus\{0\}$ such that
\[r\sum_{\xi'\sim\xi}d(\xi')\in\Z\text{ for each }\xi\in\Xi\text{ and }r\int f(x)\,dx\in \Z+\Z\alpha.\]
By (P1), $\sum_{\xi'\sim\xi}d(\xi')\neq 0$ for every $\xi\in\Xi$ and hence
$\sum_{\xi'\sim\xi}d(\xi')=n_{\xi}/r$, where $n_{\xi}\in\Z\setminus\{0\}$. It
follows that $1/r\in\sum_{i=1}^p\Q d_i$, and hence
\begin{equation}\label{nalezy}\int f(x)\,dx\in
\frac{\Z+\Z\alpha}{r}\subset\sum_{i=1}^p(\Q+\Q\alpha)d_i.
\end{equation}

Fix $\eta_1\in\Xi$ and let $\eta_2$ be an element of $\Xi$ such that
$\eta_1\nsim\eta_2$. Since $\sum_{\eta_i'\sim\eta_i}d(\eta')=n_i/r$ for
$i=1,2$,
\[\sum_{\eta_1'\sim\eta_1}n_2d(\eta_1')-\sum_{\eta_2'\sim\eta_2}n_1d(\eta_2')=0.\]
But $n_i\neq 0$ for $i=1,2$ so by (P1), there exist $\eta_1'\sim\eta_1$ and
$\eta_2'\sim\eta_2$ such that $\eta_1'-\eta_2'\in\Q+\Q\alpha$, hence
$\eta_1-\eta_2\in\Q+\Q\alpha$. Since $\eta_1$ and $\eta_2$ are arbitrary,
\[\sum_{i=1}^p d_i\xi_i\in\sum_{i=1}^p d_i(\xi_1+\Q+\Q\alpha)=S(f)\xi_1+\sum_{i=1}^p(\Q+\Q\alpha)d_i=\sum_{i=1}^p(\Q+\Q\alpha)d_i,\]
and hence
\[\sum_{i=1}^p(\xi_i+\Q+\Q\alpha) d_i=\sum_{i=1}^p(\Q+\Q\alpha) d_i.\]
Thus by (\ref{nalezy})
\[\int f(x)\,dx\in\sum_{i=1}^p(\xi_i+\Q+\Q\alpha) d_i,\]
contrary to (P2).
\end{proof}

\begin{Remark}\label{partrig}
Notice that special flows under piecewise constant functions are
partially rigid for an arbitrary irrational $\alpha$. Indeed,
suppose that $f:\T\to\R$ is a positive piecewise constant
right-continuous function for which $\Xi=\{\xi_1,\ldots,\xi_p\}$
is the set of discontinuity points of $f$ and $d_i=d(\xi_i)$ stand
for the value of the jump at $\xi_i$, $i=1,\ldots,p$. Then
\[f^{(q_n)}(x)-f^{(q_n)}(0)=-\sum_{i=1}^pd_i\sum_{j=0}^{q_n-1}\chi_{(0,x]}(\xi_i-j\alpha).\]
By the Denjoy-Koksma inequality,
\[\left|\sum_{j=0}^{q_n-1}\chi_{(0,x]}(\xi_i-j\alpha)-q_nx\right|\leq 2\]
for each $x\in\T$ and $\xi_i\in\Xi$. Since $\sum_{i=1}^pd_i=0$, it
follows that
\[f^{(q_n)}(x)-f^{(q_n)}(0)\in D:=\sum_{i=1}^pd_i\cdot\{-2,-1,0,1,2\}.\] Put $c_f=\int
f(x)\,dx$. By the Denjoy-Koksma inequality, the sequence
$(f^{(q_n)}-q_nc_f)_{n\in\N}$ is uniformly bounded, so by passing
to a subsequence, if necessary, we can assume that
\[(f^{(q_n)}-q_nc_f)_*\lambda\to P\]
in the space $\mathcal{P}(\R)$ of Borel probability measures on
$\R$, and $f^{(q_n)}(0)-q_nc_f\to\gamma$. Then $P$ is concentrated
on the finite set $\{\gamma\}+D$. Now an application of
Proposition~4.1 from~\cite{Fr-Le0} yields
\[\lambda^f\left(\left(T^f\right)_{q_nc_f}A\cap B\right)\to\int_{\R}\lambda^f\left(\left(T^f\right)_{-t}A\cap B\right)\,dP(t)\]
for every $A,B\in\mathcal{B}^f$, which implies the partial
rigidity of $T^f$ along the sequence $(q_nc_f+t_0)_{n\in\N}$,
where $t_0\in\{\gamma\}+D$ is a point for which $P(\{t_0\})>0$.
\hfill $\Box$
\end{Remark}

\section{Ratner's property}\label{sectionratner}

Let $Tx=x+\alpha$ be an irrational rotation such that $\alpha$ has
bounded partial quotients and let $f:\T\to\R^+$ be a piecewise
constant right--continuous function. In this section, we will show
that, if $f$ satisfies properties (P1) and (P2), then the special
flow $T^f$ satisfies a condition (called Ratner's property) which
emulates the Ratner condition from \cite{Ra}.

\begin{Definition} \label{defrat}(cf.\ \cite{Fr-Le,Ra,Th})
Let $(X,d)$ be a $\sigma$--compact metric space, $\mathcal{B}$ be
the $\sigma$--algebra of Borel subsets of $X$, $\mu$ a Borel
probability measure on $(X,d)$ and let $(S_t)_{t\in\R}$ be a  flow
on the space $(X,{\cal B},\mu)$. Let $P\subset\R\setminus\{0\}$ be
a finite subset and $t_0\in\R\setminus\{0\}$. The flow
$(S_t)_{t\in\R}$ is said to have {\em the property $\rat(t_0,P)$}
if for every $\vep>0$ and $N\in\N$ there exist
$\kappa=\kappa(\vep)>0$, $\delta=\delta(\vep,N)>0$ and a subset
$Z=Z(\vep,N)\in\mathcal{B}$ with $\mu(Z)>1-\vep$  such that if
$x,x'\in Z$, $x'$ is not in the orbit $x$ and $d(x,x')<\delta$,
then there are $M=M(x,x')$, $L=L(x,x')\geq N$ such that
$L/M\geq\kappa$ and there exists $\rho=\rho(x,x')\in P$ such that
\[\frac{\# \{n\in\Z\cap[M,M+L]:d(S_{nt_0}(x),S_{nt_0+\rho}(x'))<\vep\}}{L}>1-\vep.\]
Moreover, we say that $(S_t)_{t\in\R}$ has {\em the property $\rat(P)$} if the
set of $s\in\R$ such that the flow $(S_t)_{t\in\R}$ has the
$\rat(s,P)$--property is uncountable.
\end{Definition}

As it is known (see \cite{Fr-Le} or \cite{Ra} or \cite{Th}) this
Ratner's property implies that the flow $(S_t)$ is a finite
extension of any its nontrivial factor.

\begin{Lemma}\label{l1}
Assume that $f$ satisfies the property (P1). Then there exist a
finite set $F\subset\R\setminus\{0\}$, $\kappa>0$ such that for
each $N\geq1$ we can find $\delta=\delta(N)>0$ such that for each
$x\neq y$, $\|x-y\|<\delta$ there are $L=L(x,y)\geq N$,
$M=M(x,y)\geq N $ such that $L/M\geq\kappa$ and for some $\rho\in
F$ we have \beq\label{e1} f^{(n)}(x)-f^{(n)}(y)=\rho\text{ for
each integer }\,n\text{ between }M\text{ and } M+L.\eeq
\end{Lemma}
\begin{proof} Since $\alpha$ has bounded partial
quotients, there exists $0<c<1$ such that \beq\label{e2} \|j\alpha\|\geq\frac
c{|j|}\mbox{ for }j\in\Z\setminus\{0\}\;\;\text{ and
}\;\;\frac{q_{n+1}}{q_n}\leq\frac{1}{c}\text{ for all }n\in\N.\eeq Let $H$ be a
natural number such that if
$\xi_{i}-\xi_{i'}\in(\Q+\Q\alpha)\setminus(\Z+\Z\alpha)$ for some $1\leq
i,i'\leq p$ then there exists a natural number $h\leq H$ such that
\[h(\xi_{i}-\xi_{i'})=m_1+\alpha m_2\text{ with }m_1,m_2\in\Z\text{ and }|m_2|\leq H.\]
Put $\displaystyle \kappa:=\frac{c^{10}}{4pH^2(1+c^5)}$.

Given $N\geq1$ we put
$\displaystyle\delta=\frac{c^7}{2pH^2(1+c^5)N}$. Now fix $x\neq y$
such that $\|x-y\|<\delta$ and let $s\geq 1$ be a unique number so
that \beq\label{e3}\frac2{q_{s+1}}\leq\|x-y\|<\frac2{q_s}.\eeq
Then
\begin{equation}\label{jeszcze}
q_{s+1}\geq 2/\delta=\frac{4pH^2(1+c^5)}{c^7}N.
\end{equation}
We will now consider $n=q_s,q_s+1,\ldots,q_{s+4}-1$ and look at
the corresponding values $f^{(n)}(x)-f^{(n)}(y)$. We claim that
for each $q_s\leq n<q_{s+4}$ we have
\begin{equation}\label{e4}
f^{(n)}(x)-f^{(n)}(y)\in V:=\left\{\sum_{i=1}^p r_id_i:\: r_i\in
[-R,R]\cap\N\right\},
\end{equation}
where $R=2/c^5+1$.  Without loss of generality we can assume that
$x<y$. Then the function $f^{(n)}$ is piecewise constant and
\[f^{(n)}(x)-f^{(n)}(y)=\sum_{i=1}^p\sum_{j=0}^{n-1}d_i\chi_{(x,y]}(\xi_i-j\alpha).\]
Therefore we only need to show that there is a bounded number
(counted with possible multiplicities) of discontinuity points of
$f^{(q_{s+4})}$ in the interval $(x,y]$. This is however clear
because by~(\ref{e2}) and (\ref{e3}) we have
\[\|x-y\|\leq \frac{2}{c^4}\frac{1}{q_{s+4}}\] and for a
fixed $1\leq i\leq p$ and any  $0\leq j\neq j'<q_{s+4}$,
\[\|(\xi_i-j\alpha)-(\xi_i-j'\alpha)\|\geq \frac{c}{q_{s+4}},\]
hence for a fixed $1\leq i\leq p$ there are at most $R=2/c^5+1$
numbers $0\leq j<q_{s+4}$ with $\xi_i-j\alpha\in(x,y]$ and the
claim is proved.

Put $F:=V\setminus\{0\}$. Notice that it is enough to indicate a
sufficiently long integer interval contained in $[q_s,q_{s+4})$
such that for each $n$ in this interval $f^{(n)}(x)\neq
f^{(n)}(y)$. Indeed, suppose $J\subset [q_s,q_{s+4})$ a long
enough integer interval such that for each $n\in J$ we have
$f^{(n)}(x)\neq f^{(n)}(y)$. In view of (\ref{e4}),
$f^{(n)}(y)-f^{(n)}(x)\in F$ for $n\in J$. Note that if $f^{(n)}$
and $f^{(n+1)}$ have the same points of discontinuity (counted
with possible multiplicities) in the interval $(x,y]$ then
$f^{(n+1)}(y)-f^{(n+1)}(x)=f^{(n)}(y)-f^{(n)}(x)$. Since
$f^{(q_{s+4})}$ has at most $Rp$ discontinuities in $(x,y]$, we
can split $J$ into at most $Rp+1$ integer intervals so that for
every such interval $\widetilde{J}\subset J$ the sequence
$(f^{(n)}(y)-f^{(n)}(x))_{n\in \widetilde{J}}$ is constant.
Therefore if $J'$ is the longest one then $|J'|\geq|J|/(Rp+1)$ and
for some $\rho\in F$ we have $f^{(n)}(y)-f^{(n)}(x)=\rho$ for all
$n\in J'$. We will show in the sequel that there exists such an
interval $J$, with $J\subset[q_{s},q_{s+4}]$ and
$|J|\geq\min(q_{s-1}, \frac{cq_s}{2H^2})$. The bounds $M$ and
$M+L$ of an associated interval $J'$ will satisfy $M,L\geq N$ and
$L/M\geq\kappa$.

Suppose $0\leq i\leq q_{s+1}$ is the smallest so that
$f^{(q_s+i)}(x)=f^{(q_s+i)}(y)$; if there is no such $i$ we pass
directly to the final argument of this proof taking
$J=[q_s,q_{s+1}]$. We now look for the smallest $j>0$ such that
for some $1\leq k\leq p$, $\xi_k-(q_s+i+j)\alpha\in (x,y]$ (that
is we wait until a new discontinuity point comes in $(x,y]$). We
have \beq\label{e5} 0<j<q_{s+1}.\eeq Indeed, by basic properties
of the discrepancy $D_{q_{s+1}}(\alpha)$, for each $z\in\T$ in
each interval of length $2/q_{s+1}$ there must be a point of the
form $z-w\alpha$ for some $0\leq w<q_{s+1}$, and we apply this to
the interval $(x,y]$ and $z=\xi_k-(q_s+i)\alpha$ (for each
$k=1,\ldots, p$).

Hence, in view of the definition of $j$ \beq\label{e6}
f^{(q_s+i+j+1)}(x)\neq f^{(q_s+i+j+1)}(y). \eeq We now consider
$f^{(q_s+i+j+v)}(y)-f^{(q_s+i+j+v)}(x)$ for
$v=1,\ldots,q_{s+4}-(q_s+i+j)$ until (if it exists) $v_0$ for
which $f^{(q_s+i+j+v_0)}(x)=f^{(q_s+i+j+v_0)}(y)$. We have
$$ f^{(q_s+i+j+v_0)}(x)=f^{(q_s+i+j)}(x)+f^{(v_0)}(T^{q_s+i+j}x)$$
(with the same formula for $y$), so by the definition of $j$ we
obtain
$$
f^{(v_0)}(T^{q_s+i+j}x)=f^{(v_0)}(T^{q_s+i+j}y).$$ Moreover, by
the definition of $j$, we have
$\xi_k\in(T^{q_s+i+j}x,T^{q_s+i+j}y]$ (we consider intervals on
the circle) for some $1\leq k\leq p$. This implies that $f$, hence
$f^{(v_0)}$, has some discontinuity points between
$x+(q_s+i+j)\alpha$ and $y+(q_s+i+j)\alpha$. Therefore there exist
$n_i\in\N\setminus\{0\}$ for $i=1,\ldots,u$, with $1\leq u\leq p$
and some $1\leq k_1<\ldots< k_{u}\leq p$ such that
$$
0=f^{(v_0)}(T^{q_s+i+j}x)-f^{(v_0)}(T^{q_s+i+j}y)=\sum_{i=1}^{u}n_id_{k_i},$$
which by the assumption (P1) means that  there are $l\neq l'$ such that
\[\xi_{k_l}-\xi_{k_{l'}}\in(\Q+\Q\alpha)\setminus(\Z+\Z\alpha).\] Moreover,
$\|T^{q_s+i+j}x-T^{q_s+i+j}y\|<\frac2{q_s}$ and for some $0\leq w,w'< v_0$ we
have
\[\xi_{k_l}-w\alpha,\xi_{k_{l'}}-w'\alpha\in
(T^{q_s+i+j}x,T^{q_s+i+j}y],\] and in particular
\[\|(\xi_{k_l}-w\alpha)-(\xi_{k_{l'}}-w'\alpha)\|<\frac2{q_s}.\]
However
\[(\Q+\Q\alpha)\setminus(\Z+\Z\alpha)\ni\xi_{k_l}-\xi_{k_{l'}}=(m_1+
m_2\alpha)/h,\] where $m_1,m_2\in\Z$ with $|m_2|\leq H$ and $h\in\N$ with
$h\leq H$. Therefore
\begin{equation}\label{rozxi}
\begin{aligned}
\frac{\|\left(h(w-w')-m_2\right)\alpha\|}{h}&\leq\left\|\frac{m_1-(h(w-w')-m_2)\alpha}{h}\right\|\\
&=\|(\xi_{k_l}-w\alpha)-(\xi_{k_{l'}}-w'\alpha)\|<\frac{2}{q_s}.
\end{aligned}
\end{equation}
Moreover $h(w-w')-m_2\neq 0$. Indeed, if $h(w-w')-m_2=0$ then
$\|m_1/h\|<2/q_s<1/H\leq 1/h$, therefore $m_1=0$. But that gives
\[\xi_{k_l}-\xi_{k_{l'}}=\frac{m_2\alpha}{h}=(w-w')\alpha\in\Z\alpha,\]
which has been excluded. In view of~(\ref{e2}),
$\|\left(h(w-w')-m_2\right)\alpha\|\geq \frac{c}{h|w-w'|+|m_2|}$.
Hence (by (\ref{rozxi}))
$$\frac{c}{Hv_0}\leq\frac{c}{H(|w-w'|+1)}\leq\frac{c}{h|w-w'|+|m_2|}\leq \frac{2H}{q_s}$$
and finally $v_0\geq \frac{cq_s}{2H^2}$. It is now enough to take
$J=[q_s+i+j,q_s+i+j+v_0]$ to conclude.

If $f^{(q_s+i+j+v)}(y)\neq f^{(q_s+i+j+v)}(x)$ for
$v=1,\ldots,q_{s+4}-(q_s+i+j)$ then we can take
$J=[q_s+i+j,q_{s+4}]$.

Now let $J'\subset J$ be an integer interval such that
$|J'|\geq|J|/(Rp+1)$ and for some $\rho\in F$ we have
$f^{(n)}(y)-f^{(n)}(x)=\rho$ for all $n\in J'$. Let $M,L$ be
natural numbers such that $J'=[M,M+L]$. Then in view of
(\ref{jeszcze}) we obtain
\begin{eqnarray*}L&\geq&|J|/(Rp+1)\geq
\frac{c^5}{2p(c^5+1)}\min\left(q_{s+1}-q_{s},\frac{cq_s}{2H^2},q_{s+4}-2q_{s+1}-q_{s}\right)\\
&\geq
&\frac{c^5}{2p(c^5+1)}\min\left(q_{s-1},\frac{cq_s}{2H^2},q_{s+2}\right)\geq
\frac{c^7}{4pH^2(c^5+1)}q_{s+1}\geq N,\\
 M&\geq& q_{s}\geq
cq_{s+1}\geq N,
\end{eqnarray*}
and
\[\frac{L}{M}\geq\frac{c^7}{4pH^2(c^5+1)}\frac{q_{s+1}}{q_{s+4}}\geq\frac{c^{10}}{4pH^2(c^5+1)}=\kappa.\]
\end{proof}

By Corollary~\ref{slabemieszanie}, Lemma~\ref{l1} and Lemma~6 in \cite{Fr-Le}
and we have the following.

\begin{Th}\label{ratner}
Assume that $\alpha$ has bounded partial quotients and
$f:\T\to\R^+$ is a piecewise constant function satisfying
properties (P1) and (P2). Then the special flow $T^f$ has the
$\rat(\gamma,P)$--property for every $\gamma\neq 0$.
\end{Th}

Now from  Theorem~4 in \cite{Fr-Le} and the subsequent remark we have the
following.

\begin{Cor}\label{al55}
Assume that $\alpha$ has bounded partial quotients and
$f:\T\to\R^+$ is a piecewise constant function satisfying
properties (P1) and (P2). Then $T^f$ is a finite extension of each
of its non--trivial factors.
\end{Cor}

\section{Mild mixing}\label{secmm}
Assume that $T$ is an irrational rotation by $\alpha$ which has
bounded partial quotients and $f:\T\to\R^+$ is a piecewise
constant function satisfying properties (P1) and (P2).  In this
section using Ratner's property we will prove that $T^f$ is mildly
mixing.

\begin{Lemma}\label{al6}
Let $Tx=x+\alpha$ be an irrational rotation and let $f:\T\to \R$
be a positive integrable function with respect to the measure
$\lambda$. Suppose that $\ca\subset\cb^f$ is a factor of $T^f$
with finite fibers, say $k$--point in a.e.\ fiber. Assume that
$\ca$ is rigid and moreover that $T^f_{t_n}\to
E(\,\cdot\,|\mathcal{A})$, for a fixed real sequence $(t_n)$ going
to $+\infty$. Then for each $\vep>0$ we can find a family
$\{I_1,\ldots,I_l\}$ of pairwise disjoint intervals on $\T$ of the
same length $\leq\vep$ for which $\sum_{j=1}^l|I_j|>1-\vep$, and
$n_0\in\N$ such that for every $n\geq n_0$ and each $1\leq j\leq
l$ there exists a measurable set $I^{(n)}_j\subset I_j$ such that
$\lambda( I^{(n)}_j)\geq\frac12\frac1k|I_j|$ and if $x\in
I^{(n)}_j$ then there exists $m\in\N$ such that $t_n\in
f^{(m)}(x)+(-\vep,\vep)$ and $\|m\alpha\|<\vep$.
\end{Lemma}

\begin{proof}
We divide the set $\T\times(0,\vep)$ into consecutive vertical
strips (squares) $S_1,\ldots,S_l$ of width $\vep$ and the last
strip which may have smaller width. The family $I_1,\ldots,I_l$ is
defined as the base intervals of those squares. Then
$\sum_{j=1}^l|I_j|>1-\vep$. By the assumption we know that
\[\lambda^f(T^f_{-t_n}(S_j)\cap S_j)\to\int
E(S_j|\mathcal{A})\chi_{S_j}\,d\lambda^f,\] for $j=1,\ldots,l$.
Moreover
$E(S_j|\mathcal{A})(x,r)=\left(\lambda^f\right)_{\pi_{\cal
A}(x,r)}(S_j)$, where
$$
\lambda^f=\int
\left(\lambda^f\right)_{\un{x}}\,d\left(\lambda^f\right)_{\ca}(\un{x})$$
denotes the disintegration of $\lambda^f$ over the factor
$(\T^f/\ca,\left(\lambda^f\right)_{\ca})$ (and each
$\left(\lambda^f\right)_{\un{x}}$ is the uniform measure on a
$k$-element set) and $\pi_{\cal A}:\T^f\to\T^f/\mathcal{A}$ is the
natural projection. For each $S\in\mathcal{B}^f$, for a.a.\
$(x,r)\in S$ we have $\left(\lambda^f\right)_{\pi_{\cal
A}(x,r)}(S)\geq1/k$. Applied to $S=S_j$ this gives
\[\int E(S_j|{\cal A})\chi_{S_j}\,d\lambda^f\geq \frac{1}{k}\lambda^f(S_j).\]
Therefore there exists $n_0\in\N$ such that for $n\geq n_0$ we
have
\[\lambda^f(T^f_{-t_n}(S_j)\cap S_j)\geq \frac{1}{2k}\lambda^f(S_j)\] for all
$j=1,\ldots,l$. Suppose that $(x,t)\in T^f_{-t_n}(S_j)\cap S_j$
and $m\in\N$ is a unique number such that $f^{(m)}(x)\leq
t+t_n<f^{(m+1)}(x)$. Then
\[T^f_{t_n}(x,t)=(T^mx,t+t_n-f^{(m)}(x))\in  S_j\subset X\times(0,\vep)\]
and hence $t_n\in f^{(m)}(x)+(-\vep,\vep)$. Since
$(x,t),(T^mx,t+t_n-f^{(m)}(x))\in S_j$ and the width of $S_j$ is
less  than $\vep$, we obtain $\|m\alpha\|=\|T^mx-x\|<\vep$. Taking
$ I_j^{(n)}=\{x\in I_j:\exists_{0\leq t<\vep}\,(x,t)\in
T^f_{-t_n}(S_j)\cap S_j\}$ completes the proof.
\end{proof}

Assume that $\alpha$ has bounded partial quotients and
$f:\T\to\R^+$ is a piecewise constant function satisfying
properties (P1) and (P2).

\begin{Lemma}\label{al7}
Under the assumptions on $\mathcal{A}$ of Lemma~\ref{al6} there
exist a finite set $0\notin F\subset\R$ and  $0<\xi<1$ such that
for each (sufficiently small) $\vep>0$ and $n\geq n_0$,
($n_0=n_0(\vep)$) there exists a family $(I_j)_{1\leq j\leq l}$ of
pairwise disjoint intervals of length $\leq\vep$ such that
$\sum_{j=1}^l|I_j|>1-\vep$ and there exists a family
$\left(\widetilde{I}_j^{(n)}\right)_{1\leq j\leq l}$ of measurable
sets such that $\widetilde{I}^{(n)}_j\subset I_j$  and
$\lambda(\widetilde{I}^{(n)}_j)>\xi|I_j|$ for $j=1,\ldots,l$, and
for every $x\in \widetilde{I}^{(n)}_j$ there exists $s\in\N$ for
which
\[f^{(s)}(x)\in t_n+F+(-\vep,\vep)\text{ and }\|s\alpha\|<\vep.\]
\end{Lemma}

Let us insist on the fact that the novelty of Lemma~\ref{al7} after
Lemma~\ref{al6} is the fact that $0$ does not belong to $F$. Roughly speaking,
the  proof of Lemma~\ref{al7} goes along the following lines: a little $\vep$
and a large $n$ are fixed. Since $f$ is piecewise constant, by Lemma~\ref{al6}
and its proof, there is a maximal family of intervals $ K_1, K_2,\ldots,K_l $
such that, for each $j$, there exists $s$ such that for all $ x\in K_j $, $ |
f^{(s)}(x)-t_n | < \vep$. The measure of $ \cup_{j=1}^L K_j $ is not too small.
Also, the distances between these $K_j$ are not too small. We  associate to
each $K_j$ a right adjacent interval $K'_j$. These new intervals will do the
job.

\begin{proof}
Let $a$ and $b$ stand for the minimal  and maximal value of $f$
respectively. Let $F$ be the set of sums
 \[d_{k_1}+d_{k_2}+\ldots+d_{k_u}\]
such that $1\leq k_1<\ldots<k_u\leq p$ satisfy
$\xi_{k_i}-\xi_{k_{i'}}\notin(\Q+\Q\alpha)\setminus(\Z+\Z\alpha)$ for all
$1\leq i,i'\leq u$. By (P1), $0\notin F$. Let
\[F_1=\{n_0a+n_1d_1+\ldots+n_pd_p:0<|n_0|\leq 2pb/a, |n_i|\leq 2pb/a+1,i=1,\ldots,p\}.\]
In view of (\ref{wznien}), $0\notin F_1$. Recall that (see
Remark~\ref{nieciagl}) for every $m\in\N$ the set
$\{\xi_i-j\alpha:1\leq i\leq p,0\leq j<m\}$ is the set of
discontinuities of $f^{(m)}$.

From Lemmas~\ref{al1} and \ref{al2} there exist $c_1,c_2>0$ such that for any
natural $m$ the lengths of subintervals given by the partition of $\T$ into the
points of discontinuity of $f^{(m)}$ are at most $c_2/m$ and if
$\xi_i-\xi_j\in\Q+\Q\alpha$ and $\xi_i-r\alpha\neq\xi_j-s\alpha$ then
\begin{equation}\label{szadl}
\|(\xi_i-r\alpha)-(\xi_j-s\alpha)\|\geq
\frac{c_1}{m}\text{ for all }0\leq r,s<m.
\end{equation}

Fix $0<\vep<\min(a/2,1/2)$ such that
\begin{equation}
(-2\vep,2\vep)\cap (F\cup F_1)=\emptyset.
\end{equation}
By Lemma~\ref{al6} we get the family of intervals
$\{I_1,\ldots,I_l\}$ and the natural number $n_0=n_0(\vep)$. We
consider now $n\geq n_0$ such that
\begin{equation}\label{duzn}
t_n>\max(16klc_2b,3\vep).
\end{equation}
 We now estimate for how many $j$'s the graph of $f^{(j)}$ has a non-empty intersection
with the rectangle $\T\times[t_n-\vep, t_n+\vep]$. If $m$ is the
maximal, and $m_1$ is the minimal such a $j$, then $f^{(m)}\geq
am$ and $f^{(m_1)}\leq m_1b$, so
\begin{equation}\label{do1}
m_1\geq m\cdot\frac{a}{2b}\;\;\;\;\text{ and }\;\;\;\;t_n/2\leq m_1b.
\end{equation}
Furthermore,  the length of any interval on which $f^{(s)}$ is
constant for some $m_1\leq s\leq m$ is at most $c_2/m_1$.
Moreover, for every $m_1\leq s\leq m$
\begin{equation}\label{text}
\begin{aligned}
&\text{ if }0<x'-x<c_1/m,\;\;x,x'\text{ are continuity points of }f^{(s)}\\&
\text{ and }f^{(s)}\text{ is not constant on }[x,x']\text{ then }
f^{(s)}(x')-f^{(s)}(x)\in F.
\end{aligned}
\end{equation}
Indeed,
\[f^{(s)}(x)-f^{(s)}(x')=\sum_{i=1}^p\sum_{j=0}^{s-1}d_i\chi_{(x,x')}(\xi_i-j\alpha)\]
and by (\ref{szadl}) if $\xi_i-j\alpha,\xi_{i'}-j'\alpha\in(x,x')$ then either
$\xi_i-\xi_{i'}\notin\Q+\Q\alpha$ or $\xi_i-j\alpha=\xi_{i'}-j'\alpha$. It
follows that for every $1\leq i\leq p$ there exists at most one $0\leq j<s$
such that $\xi_i-j\alpha\in(x,x')$, and hence
\[f^{(s)}(x)-f^{(s)}(x')=d_{k_1}+d_{k_2}+\ldots+d_{k_u},\]
where $1\leq k_1<\ldots<k_u\leq p$ and
$\xi_{k_i}-\xi_{k_{i'}}\notin(\Q+\Q\alpha)\setminus(\Z+\Z\alpha)$
for all $1\leq i,i'\leq u$.

Let us consider the set
\[\widehat{I}^{(n)}_j:=\bigcup_{m_1\leq s\leq m}\wne
\{x\in I_{j}: f^{(s)}(x)\in (t_n-\vep,t_n+\vep),
\|s\alpha\|<\vep\}.\] Then $\widehat{I}^{(n)}_j$ is the union of
disjoint consecutive intervals $J_1,\ldots,J_{N+1}$, i.e.\ $J_i$
lies on the left hand side of $J_{i+1}$ for all $1\leq i\leq N$.
Since $ I^{(n)}_j\subset \bigcup_{i=1}^{N+1}\overline{J_i}$, we
have $\lambda(\widehat{I}^{(n)}_j)>\frac{1}{2k}|I_j|$.

Suppose that $x\in \widehat{I}^{(n)}_j$ and
$f^{(s)}(x),f^{(s')}(x)\in (t_n-\vep,t_n+\vep)$ for some $m_1\leq
s\leq s'\leq m$. Then
\[(s'-s)a\leq f^{(s'-s)}(T^{s}x)=f^{(s')}(x)-f^{(s)}(x)\in[0,2\vep)\subset[0,a),\]
hence $s=s'$.

Suppose that $x\in J_i$ and $s$ is a natural number such that
$f^{(s)}(x)\in t_n+(-\vep,\vep)$. Then $f^{(s)}$ is constant on
$J_i$. Indeed, otherwise there exist continuity points $x,x'\in
J_i$ of $f^{(s)}$ such that $0<x'-x<c_1/m$ and $f^{(s)}$ is not
constant on $[x,x']$; then by (\ref{text}),
\[ f^{(s)}(x')-f^{(s)}(x)\in F\cap(-2\vep,2\vep),\]
which is impossible. It follows also that for every $i=1,\ldots,N+1$ there
exists $m_1\leq s_i\leq m$ such that $\|s_i\alpha\|<\vep$ and $f^{(s_i)}(x)\in
t_n+(-\vep,\vep)$ for all $x\in J_i$.  It is also clear that the ends of the
interval $J_i$ are discontinuities of $f^{(s_i)}$ for $i=1,\ldots,N+1$, maybe
except for the left hand end of $J_1$ and the right hand end of $J_{N+1}$.
Therefore $|J_i|\leq c_2/m_1$ for $i=1,\ldots,N+1$. Moreover,
\begin{equation}\label{dluganie}
\frac{1}{4kl}\leq\frac{1}{2k}|I_j|\leq\lambda(I^{(n)}_j)\leq
\lambda(\bigcup_{i=1}^{N+1}J_i)\leq (N+1)\frac{c_2}{m_1}.
\end{equation}
Hence, by (\ref{do1}) and (\ref{duzn}) we obtain
\[N+1\geq\frac{m_1}{4klc_2}\geq\frac{t_n}{8klc_2b}\geq 2\]
and hence $N\geq 1$.

Notice that for every  $1\leq i\leq N$ the distance between $J_i$ and $J_{i+1}$
is at least $c_1/m$. Otherwise, there exist continuity points $x\in J_i$,
$x'\in J_{i+1}$ for $f^{(s_{i+1})}$ such that $0<x'-x<c_1/m$. Since
$f^{(s_{i+1})}$ is not constant on $[x,x']$, by (\ref{text}), we have
$f^{(s_{i+1})}(x')-f^{(s_{i+1})}(x)\in F$. Then
\begin{eqnarray*}
f^{(s_{i+1}-s_{i})}(T^{s_{i}}x)&=&f^{(s_{i+1})}(x)-f^{(s_{i})}(x)\in
f^{(s_{i+1})}(x')-f^{(s_{i})}(x)-F\\&\subset&
-F+(-2\vep,2\vep)\subset(-2pb,2pb),
\end{eqnarray*}
since $F\subset(-pb,pb)$. It follows that $|s_{i+1}-s_{i}|\leq 2pb/a$ and hence
\[f^{(s_{i+1}-s_{i})}(T^{s_{i}}x)=(s_{i+1}-s_{i})a+n_1d_1+\ldots+n_p d_p,\]
where $|n_j|\leq 2pb/a$ for $j=1,\ldots,p$. If $s_{i+1}\neq s_i$
then
\[(-2\vep,2\vep)\ni f^{(s_{i+1})}(x')-f^{(s_{i})}(x)\in f^{(s_{i+1}-s_{i})}(T^{s_{i}}x)+F\subset F_1,\]
which is a contradiction. Therefore $s_{i+1}=s_{i}$ and then
\[f^{(s_{i})}(x')-f^{(s_{i})}(x)\in F\cap(-2\vep,2\vep),\]
which is also impossible.

 For every
$i=1,\ldots,N$ let $J'_i$ stand for the open  right adjacent to
$J_i$ interval of length $c_1/m$. If $x'\in J'_i$ is a continuity
point for $f^{(s_i)}$ then we can find $x\in J_i$ such that
$0<x'-x<c_1/m$. So we have
\begin{equation}\label{istota}
f^{(s_i)}(x')\in f^{(s_i)}(x)+F\subset F+t_n+(-\vep,\vep),
\end{equation}
by (\ref{text}), because $f^{(s_i)}$ is not constant on $[x,x']$. Clearly,
(\ref{istota}) holds also for possible discontinuity points of $f^{(s_i)}$.
Moreover
\[\lambda(\bigcup_{i=1}^NJ_i')\geq \frac{N}{N+1}\frac{c_1/m}{c_2/m_1}\lambda(\bigcup_{i=1}^{N+1}J_i)\geq
\frac{c_1a}{8kc_2b}|I_j|,\] by (\ref{dluganie}) and (\ref{do1}). Set
\[\widetilde{I}^{(n)}_j:=\bigcup_{i=1}^NJ_i'.\]
Then $\lambda(\widetilde{I}^{(n)}_j)>\xi|I_j|$, where
$\xi=\frac{c_1a}{8kc_2b}$.
\end{proof}

\begin{Lemma}\label{al8}
Under the assumptions on $\mathcal{A}$ of Lemma~\ref{al6} there
exist a finite set $F$, with $0\notin F\subset\R$ and a real
$\xi$, with $0<\xi\leq 1$ such that for each $\vep>0$ there exist
$n_0\in\N$ and a finite family $\{P_{i,j,\vep}\}_{i,j}$ of
pairwise disjoint rectangles in $\T^f$ with diameter less than
$\vep$ and whose union is of $\lambda^f$-measure at least $1-\vep$
such that for each $P_{i,j,\vep}$ and $n\geq n_0$ there exists
$\widetilde{P}^{(n)}_{i,j,\vep}\subset P_{i,j,\vep}$ with
$\lambda^f(\widetilde{P}^{(n)}_{i,j,\vep})\geq\xi\lambda^f(P_{i,j,\vep})$
such that if $(x,t)\in \widetilde{P}^{(n)}_{i,j,\vep}$ then
\[d(T^f_{t_n+\rho}(x,t),(x,t))<\vep\text{ for some }\rho\in F.\]
\end{Lemma}
\begin{proof}
We fix $\vep>0$ and consider the family $\{I_i\}_{i=1}^l$ as in
Lemma~\ref{al6} (taking $\vep/2$ instead of $\vep$). Let us
consider the family of all rectangles of the form
$P_{i,j,\vep}=I_i\times (j\vep/2,(j+1)\vep/2)$ which are subsets
of\begin{eqnarray*}\T^f_\vep&:=&\{(x,t)\in\T^f:\vep<t<f(x)-\vep\}\\&&\setminus
\left(\bigcup_{u=1}^p[\xi_u-\vep,\xi_u+\vep]\times[\liminf_{x\to\xi_u}f(x)-\vep,\limsup_{x\to\xi_u}f(x)]\right).
\end{eqnarray*}
Taking $\vep$ small enough, we can make $\lambda^f\left(\cup_{i,j}
P_{i,j,\vep}\right)$ arbitrarily close to $1$. By Lemma~\ref{al7},
there exist $0<\xi<1$ (independent of $\vep$) and $n_0$ such that
for each $n\geq n_0$ there exists a family of measurable sets
$\{\widetilde{I}^{(n)}_1,\ldots,\widetilde{I}^{(n)}_l\}$  such
that $\widetilde{I}^{(n)}_i\subset I_i$  and
$\lambda(\widetilde{I}^{(n)}_i)>\xi|I_i|$ for $i=1,\ldots,l$ and
for every $x\in \widetilde{I}^{(n)}_i$ there exist $s\in\N$ and
$\rho\in F$ for which
\begin{equation}\label{epspol}
f^{(s)}(x)\in t_n+\rho+(-\vep/2,\vep/2)\text{ and
}\|s\alpha\|<\vep/2.
\end{equation}
Suppose that $(x,t)\in P_{i,j,\vep}$ and $x\in
\widetilde{I}_i^{(n)}$. Let $s$ be a natural numer for which
(\ref{epspol}) holds. Since $(x,t)\in \T^f_\vep$ and
$\|s\alpha\|<\vep$, we have $(T^sx,t)\in\T^f$ and
$\vep<t<f(T^sx)-\vep$. Let $\vep_n=t_n+\rho-f^{(s)}(x)$. Then
$|\vep_n|<\vep/2$ and $0<t+\vep_n<f(T^sx)$. Therefore
$(T^sx,t+\vep_n)\in\T^f$ and
\[(T^sx,t+\vep_n)=T^f_{\vep_n+f^{(s)}(x)}(x,t)=T^f_{t_n+\rho}(x,t).\]
Consequently
\[d(T^f_{t_n+\rho}(x,t),(x,t))=\|T^sx-x\|+|\vep_n|=\|s\alpha\|+|\vep_n|<\vep.\]
Now taking $\widetilde{P}^{(n)}_{i,j,\vep}=\{(x,t)\in
P_{i,j,\vep}:x\in\widetilde{I}_i^{(n)}\}$ completes the proof.
\end{proof}

\begin{Lemma}\label{al9}
Under the assumptions on $\mathcal{A}$ of Lemma~\ref{al6} there
exist a finite set $F$, with $0\notin F\subset\R$ and a real
$\xi$, with $0<\xi\leq 1$ such that for every $A\in\cb^f$ we have
\[\sum_{\rho\in F}\lambda^f(T^f_{-t_n-\rho}A\cap
A)\geq \frac{\xi}{2}\lambda^f(A)\] for each $n$ large enough.
\end{Lemma}
\begin{proof}
First suppose that $A\subset \T^f$ is a rectangle. For every $\vep>0$ let
\[A_{\vep}=\bigcup_{i,j}\{P_{i,j,\vep}\subset A:d(P_{i,j,\vep},\partial A)>\vep\}.\]
Fix $\vep>0$ such that
$\lambda^f(A_{\vep})\geq\frac{1}{2}\lambda^f(A)$. By
Lemma~\ref{al8}, for all $n$ large enough there exists
$\widetilde{P}^{(n)}_{i,j,\vep}\subset P_{i,j,\vep}$ with
$\lambda^f(\widetilde{P}^{(n)}_{i,j,\vep})\geq\xi\lambda^f(P_{i,j,\vep})$
such that if $(x,t)\in \widetilde{P}^{(n)}_{i,j,\vep}$ then
\[d(T^f_{t_n+\rho}(x,t),(x,t))<\vep\text{ for some }\rho\in F.\]
Let
\[A_{\vep,n}=\bigcup_{i,j}\{\widetilde{P}^{(n)}_{i,j,\vep}\subset A_{\vep}:d(\widetilde{P}^{(n)}_{i,j,\vep},
\partial A)>\vep\}.\]
Then $\lambda^f(A_{\vep,n})>\xi\lambda^f(A_{\vep})$. If for each
$\rho\in F$ we put
\[A_{\vep,n}^\rho=\{(x,t)\in A_{\vep,n}:d(T^f_{t_n+\rho}(x,t),(x,t))<\vep\}\]
then $A_{\vep,n}^{\rho}\subset T^f_{-(t_n+\rho)}A\cap A$ and
$A_{\vep,n}=\bigcup_{\rho\in F}A_{\vep,n}^\rho$, which yields
\begin{equation}\label{najoi}
 \frac{\xi}{2}\lambda^f(A)\leq
\lambda^f(A_{\vep,n})\leq\sum_{\rho\in
F}\lambda^f(A_{\vep,n}^\rho)\leq\sum_{\rho\in
F}\lambda^f(T^f_{-t_n-\rho}A\cap A).
\end{equation}
From this inequality, it is immediate that
$$ \frac{\xi}{2}\lambda^f(A)\leq
\sum_{\rho\in F}\lambda^f(T^f_{-t_n-\rho}A\cap A).
$$
is still true for any finite union of pairwise disjoint
rectangles. By standard approximation of a Borel set in
$\T\times\R$ by such  unions, the inequality extends to any
element of $\cb^f$.
\end{proof}

\begin{prooftwgl}
Suppose, contrary to our claim, that $T^f$ is not mildly mixing. By
Lemma~\ref{al5}, there exists a nontrivial rigid factor
$\mathcal{A}\subset\mathcal{B}^f$ of $T^f$ and $(t_n)_{n\in\N}$ with
$t_n\to+\infty$ such that $T^f_{t_n}\to E(\,\cdot\,|\mathcal{A})$ weakly. From
Corollary~\ref{al55} it follows that $\ca\subset\cb^f$ is a factor of $T^f$
with finite fibers. Using Lemma~\ref{al9}, for every $A\in\ca$ we have
\[\sum_{\rho\in F}\lambda^f(T^f_{-t_n-\rho}A\cap A)\geq
\frac{\xi}{2}\lambda^f(A)\] for each $n$ large enough. Since
$\lambda^f(T^f_{t_n}A\triangle A)\to 0$, we conclude that
\begin{equation}\label{wzkon}
\sum_{\rho\in F}\lambda^f(T^f_{-\rho}A\cap A)\geq
\frac{\xi}{2}\lambda^f(A).
\end{equation} for every $A\in\ca$. Let us
consider the flow $T^f$ on $\T^f/\ca$. First notice that
$\T^f/\ca$ cannot be finite. Otherwise, since the extension
$\T^f\to \T^f/\ca$ is finite, $\T^f$ is finite, a contradiction.
As $T^f$ on $\T^f/\ca$ is ergodic and $\T^f/\ca$  is not finite, a
Rokhlin lemma for flows insures the existence of $A\in\ca$ of
positive measure such that $T^f_{-\rho}A\cap A=\emptyset$ for all
$\rho\in F$. \hfill $\Box$
\end{prooftwgl}
{\em E-mail addresses:} \parbox[t]{8cm}{
 fraczek@mat.uni.torun.pl\\
mlem@mat.uni.torun.pl\\
lesigne@univ-tours.fr }


\footnotesize

\begin{thebibliography}{99}
\bibitem{Aar}J.\ Aaronson, {\em An introduction to infinite ergodic
theory}, Mathematical Surveys and Monographs, {\bf 50} AMS,
Providence, RI, 1997.
\bibitem{Arn}V.I.\ Arnold, {\em Topological and ergodic properties
of closed $1$-forms with incommensurable periods}, (Russian)
Funktsional.\ Anal.\ i Prilozhen.\ {\bf 25} (1991), 1-12;
translation in Funct.\ Anal.\ Appl.\ 25 (1991), 81-90.
\bibitem{Co-Fo-Si}I.P.\ Cornfeld, S.V.\ Fomin, Ya.G.\ Sinai, {\em Ergodic Theory},
Springer-Verlag, New York, 1982.
\bibitem{Fr-Le0}K.\ Fr\k{a}czek, M.\ Lema\'nczyk, {\em A class of special
flows over irrational rotations which is disjoint from mixing flows}, Ergod.\
Th.\ Dynam.\ Sys.\ {\bf 24} (2004), 1083-1095.
\bibitem{Fr-Le}K.\ Fr\k{a}czek, M.\ Lema\'nczyk, {\em
On mild mixing of special flows over irrational rotations under
piecewise smooth functions}, Ergod.\ Th.\ Dynam.\ Sys.\ {\bf 26}
(2006), 719-738.
\bibitem{Fr-Le1}K.\ Fr\k{a}czek, M.\ Lema\'nczyk, {\em
Smooth singular flows in dimension 2 with the minimal self-joining
property}, preprint.
\bibitem{Fur}H.\ Furstenberg, {\em Recurrence in ergodic theory
and combinatorial number theory}, Princeton University Press, Princeton, N.J.,
1981.
\bibitem{Fu-We}H.\ Furstenberg, B.\ Weiss, {\em The finite multipliers of
infinite ergodic transformations. The structure of attractors in dynamical
systems} (Proc.\ Conf., North Dakota State Univ., Fargo, N.D., 1977), Lecture
Notes in Math.\ {\bf 668}, Springer, Berlin, 1978, 127--132.
\bibitem{Ga-Le-Li}P.\ Gabriel, M.\ Lema\'nczyk,
P.\ Liardet, {\em Ensemble d'invariants pour les produits croisés
de Anzai},  Mém.\ Soc.\ Math.\ France (N.S.) {\bf 47} (1991), 102
pp.
\bibitem{Gl}E.\ Glasner, {\em Ergodic theory via joinings}, Mathematical Surveys and Monographs,
{\bf 101} AMS, Providence, RI, 2003.
\bibitem{Her}M.\ Herman, {\em Sur la conjugaison
diff\'erentiable des diff\'eomorphismes du cercle \`a des
rotations}, Publ.\ Mat.\ IHES {49} (1979), 5-234.
\bibitem{KaR}A.\ Katok, {\em Cocycles, cohomology and combinatorial
constructions in ergodic theory} (in collaboration with E. A. Robinson, Jr.),
in Smooth Ergodic Theory and its applications, Proc.\ Symp.\ Pure Math., {\bf
69} (2001), 107-173.
\bibitem{Ko}A.V.\ Kochergin, {\em On the absence of mixing in special
flows over the rotation of a circle and in flows on a
two-dimensional torus}, Dokl.\ Akad.\ Nauk SSSR {\bf 205} (1972),
949-952.
\bibitem{Ch}Y.\ Khinchin, {\em Continued Fractions}, Chicago Univ.\ Press, Chicago, 1964.
\bibitem{Le-Pa}M.\ Lema\'nczyk, F.\ Parreau, {\em Lifting mixing
properties by Rokhlin cocycles}, preprint.
\bibitem{vNe}J.\  von Neumann, {\em Zur Operatorenmethode in der klassischen Mechanik},
Ann.\ of Math.\ (2) {\bf 33} (1932), 587-642.
\bibitem{Ra}M.\ Ratner, {\em
Horocycle flows, joinings and rigidity of products}, Ann.\ of
Math.\ (2) {\bf 118} (1983), 277-313.
\bibitem{Sch}K.\ Schmidt, {\em Asymptotic properties of unitary representations and
mixing},  Proc.\ London Math.\ Soc.\ (3) {\bf 48} (1984), 445-460.
\bibitem{Sch-Wa}K.\ Schmidt, P.\ Walters, {\em Mildly mixing actions of
locally compact groups}, Proc.\ London Math.\ Soc.\ (3) {\bf 45}
(1982), 506-518.
\bibitem{Th}J.-P.\ Thouvenot, {\em Some properties and applications
of joinings in ergodic theory}, in Ergodic theory and its
connections with harmonic analysis (Alexandria, 1993), London
Math.\ Soc.\ Lecture Note Ser., 205, Cambridge Univ.\ Press,
Cambridge, 1995, 207-235.
\end{thebibliography}
\normalsize
\end{document}